\documentclass[11pt]{elsarticle}
\usepackage{amssymb}
\usepackage{theorem}
\usepackage[british]{babel}
\usepackage[mathscr]{euscript}
\usepackage{graphicx}
\usepackage{amsmath}
\usepackage{array}
\usepackage{subfigure}
\usepackage{xspace}
\usepackage{wrapfig}

\usepackage[left=2.3cm,top=2cm,right=2.3cm,nohead,foot=1cm]{geometry}
\newcommand{\rmd}{{\mathrm{d}}}
\newcommand{\Real}{\mathbb{R}}
\newcommand{\Pbb}{\mathbb{P}}

\usepackage{amsfonts}
\usepackage{amsmath}
\usepackage[all]{xy}
\usepackage{graphics}
\usepackage{amssymb}
\usepackage{url}
\theoremstyle{plain}
\newtheorem{Theorem}{Theorem}[section]
\newtheorem{Corollary}{Corollary}[section]
\newtheorem{Lemma}{Lemma}[section]
\newtheorem{Proposition}{Proposition}[section]

\newtheorem{remark}{Remark}
\newtheorem{Example}{Example}

\begin{document}

\begin{frontmatter}

\title{Estimating level sets of a distribution function using  a plug-in method: \\ a multidimensional extension}

\author{Elena Di Bernardino\footnote{Universit\'{e} de Paris X Ouest, UFR SEGMI, 200 Avenue de la R\'{e}publique, 92000 Nanterre,  France.   elena.di-bernardino@univ-lyon1.fr,  http://isfaserveur.univ-lyon1.fr/\~elena.dibernardino/.},  Thomas Lalo\"{e}\footnote{Universit\'{e} de Nice Sophia-Antipolis, Laboratoire J-A Dieudonn\'{e}, Parc Valrose, 06108 Nice Cedex 02.  Thomas.Laloe@unice.fr, http://math.unice.fr/\~laloe/}}

\begin{abstract}
\begin{sloppypar}
This paper deals with the problem of estimating the level sets ${L(c)= \{F(x) \geq c \}}$, with $c \in (0,1)$, of an unknown distribution function $F$ on $\mathbb{R}^d_+$. A plug-in approach is followed. That is, given a consistent 
estimator $F_n$ of $F$, we estimate $L(c)$ by $L_n(c)= \{F_n(x) \geq c \}$.  We state consistency  results with respect to the Hausdorff distance and the volume of the symmetric difference. These results  can be considered as 
generalizations of results  previously obtained,  in a bivariate framework, in Di Bernardino \emph{et al.} (2011\nocite{DiBernardino}).  Finally we investigate the effects  of scaling  data on our consistency results. 
\end{sloppypar}
\end{abstract}

\begin{keyword}
  Level sets, multidimensional distribution function, plug-in estimation, Hausdorff distance.
\end{keyword}
\end{frontmatter}

\section*{Introduction}\label{Introduction}

In this present paper, we consider the problem of estimating the level sets of a $d$-variate distribution function. To this aim, we generalize the results obtain in a previous paper (Di Bernardino \emph{et al.},
2011\nocite{DiBernardino}).\\

As yet remarked  in Di Bernardino \emph{et al.} (2011\nocite{DiBernardino}), considering the level sets of a distribution function,  the commonly assumed property of compactness for these sets is no more reasonable.
Then, differing  from the classical literature  (Cavalier, 1997\nocite{Cavalier};  Cuevas and Fraiman, 1997\nocite{Cuevas1};  Ba\'{\i}llo \textit{et al.}, 2001\nocite{Baillo}; Ba\'{\i}llo, 2003\nocite{Baillo2};
Cuevas \textit{et al.}, 2006\nocite{Cuevas2}; Biau \textit{et al.}, 2007\nocite{Biau}; Lalo\"{e}, 2009\nocite{Laloe}),  we  need to work in a non-compact setting and this requires special attention in the statement of our problem.\\

We follow the same general approach than in Di Bernardino \emph{et al.} (2011\nocite{DiBernardino}), and we will keep as much as possible the same notation. Considering a consistent estimator $F_n$ of the distribution function
$F$, we propose a plug-in approach (e.g. see Ba\'{\i}llo \textit{et al.}, 2011\nocite{Baillo};  Rigollet and Vert\nocite{Rigollet}, 2009; Cuevas \textit{et al.}, 2006\nocite{Cuevas2}) to estimate the level set

$$L(c)= \{x \in \Real^d_+: F(x) \geq c \}, \quad \mbox{ for } c \in (0,1),$$
by
$$L_n(c)= \{x \in \Real^d_+:  F_n(x) \geq c \}, \quad \mbox{ for } c \in (0,1).$$

\vspace{0.12cm}
The  regularity  properties of $F$ and $F_n$ as well as the consistency properties of $F_n$ will be specified in the statements of our theorems.\\

As in  Di Bernardino \emph{et al.} (2011)\nocite{DiBernardino} our consistency results are stated with respect to   two criteria of ``physical proximity'' between sets: the Hausdorff  distance and  volume of the symmetric
difference.  If the consistency in term of the Hausdorff distance is a trivial extension of Theorem 2.1 in Di Bernardino \emph{et al.} (2011\nocite{DiBernardino}) (see Theorem \ref{convergenza su insiemi cresenti} below), things
are a little more complex for the  volume of the symmetric difference. In particular, in this latter case the convergence rate  suffers from the well-known \emph{curse of dimensionality}  (see Theorem  \ref{convergenza lambda su insiemi cresenti}).\\

A second aim  of this paper is to analyze  the effects of scaling data on our consistency results (see Theorem \ref{meme A}).\\

The paper is organized as follows.  We introduce some notation, tools and technical assumptions in Section \ref{Notations}.  Consistency  and  asymptotic properties of our estimator of $L(c)$ are given in Sections
\ref{Hausdorff distance} and \ref{lambda consistency}. Section \ref{variation echelle} is devoted to investigate the effects  of scaling  data on our consistency results.   Finally, proofs are postponed to Section \ref{proofs}.\\

\section{Notation and preliminaries}\label{Notations}

In this section we introduce some notation and tools which will be useful later.\\

\begin{sloppypar}
Let $\mathbb{N}^*=\mathbb{N}\setminus\{0\}$, $\Real^*_+ = \Real_+ \setminus\{0\}$ and  ${\mathbb{R}^{d}_{+}}^{*} = \mathbb{R}^{d}_{+} \setminus\{0\}$. Let $\mathcal{F}$ be  the set of continuous distribution functions
${\mathbb{R}^d_+ \to [0,1]}$ and $\textbf{X}:=(X_1,X_2,\hdots,X_d)$  a random vector with  distribution function $F\in \mathcal{F}$.   Given an $i.i.d$ sample $\{\textbf{X}_i\}_{i=1}^{n}$ in $\mathbb{R}^{d}_{+}$ with distribution
function $F \in \mathcal{F}$, we denote by $F_n$ an estimator of $F$  based on this finite sample.  \end{sloppypar}

\vspace{0.2cm}
Define, for $c \in (0,1)$, the \textit{upper $c$-level set} of $F \in \mathcal{F}$ and its plug-in estimator
$$L(c) = \{x\in \mathbb{R}^{d}_{+}: F(x) \geq c \},  \quad L_n(c) = \{x\in \mathbb{R}^{d}_{+}: F_n(x) \geq c \},$$

and
$$\{F=c\}= \{x\in \mathbb{R}^{d}_{+}: F(x) = c\}.$$

In addition, given $T>0$, we set
$$ L(c)^T = \{x\in [0,T]^d : F(x) \geq c \},\,\,\, L_n(c)^T = \{x \in [0,T]^d : F_n(x) \geq c \},$$
$$\{F=c\}^T = \{x \in [0,T]^d : F(x) = c\}.$$

Given a set $A\subset\mathbb{R}^{d}_{+}$  we denote by $\partial A$ its boundary, and by $\beta\,A$ the scaled set $\{\beta\, x, \,\mbox{ with }\, x \in A\}$.\\

Note that,  in the presence of a plateau at level $c$,  $\{F=c\}$ can be a portion of quadrant $\mathbb{R}^{d}_{+}$ instead of a set of Lebesgue measure null in $\mathbb{R}^{d}_{+}$. In the following  we introduce suitable conditions in order to avoid this situation.\\

\begin{sloppypar}
We denote by $B(x, \rho)$ the closed ball centered on $x \in \mathbb{R}^{d}_{+}$ and with positive radius $\rho$. Let  ${B(S, \rho) = \bigcup_{x \in S} B(x, \rho)}$, with  $S$ a closed set of $\mathbb{R}^d_+$.
\end{sloppypar}
For $r>0$ and  $\zeta >0$, define  $$E = B(\{x \in \mathbb{R}^{d}_{+} :\,  \mid F(x)-c\mid \leq r\}, \, \zeta),$$
and, for a twice differentiable function $F$,
\begin{equation*}\label{m_delta}
m^{\triangledown} =\inf_{x \in E} \|(\nabla F)_x\|, \qquad  M_H = \sup_{x \in E} \|(H F)_{x}\|,
\end{equation*}
where $(\nabla F)_x$ is the gradient vector of $F$ evaluated at $x$ and  $ \|(\nabla F)_x\|$ its Euclidean norm,  $(H F)_{x}$ the Hessian matrix evaluated in $x$ and $\|(H F)_{x}\|$ its matrix norm induced by the Euclidean norm.\\

For sake of completeness, we recall that if $A_1$ and $A_2$ are compacts sets in $\mathbb{R}^d_+$, the Hausdorff  distance
between $A_1$ and $A_2$ is defined by

$$d_H(A_1,A_2) = \max\left\{\sup_{x \in A_1} d(x,A_2) , \sup_{x \in A_2} d(x,A_1) \right\},$$
where $d(x,A_2)= \inf_{y \in A_2} \parallel x-y \parallel$.\\

The above expression is well defined even when $A_1$ and $A_2$ are just closed (not necessarily compacts) sets but, in this case, the value $d_H(A_1,A_2)$ could be infinity. Then in our setting, in order to avoid these situations,
we introduce the following assumption.\\

\vspace{0.1cm}
\begin{itemize}
\label{H}
\item[$\mathbf{H}$:]   There exist  $\gamma > 0$ and $A > 0$ such that, if $\,|\,t - c\,|\, \leq \gamma$ then $\forall \, \,$ $T > 0\,$ such that  $\{F=c\}^T  \neq   \emptyset$ and   $\{F=t\}^T  \neq   \emptyset, \,\,$
$$d_H (\{F=c\}^T, \{F=t\}^T) \leq \,  A \, \,|\,t - c\,|\,.$$
\end{itemize}

\vspace{0.1cm}
For further details about this assumption the interest reader is referred to  Di Bernardino \emph{et al.} (2011)\nocite{DiBernardino}, Cuevas \textit{et al.} (2006)\nocite{Cuevas2}, Tsybakov \nocite{Tsybakov}(1997). Remark that  a sufficient condition for Assumption $\mathbf{H}$ can be obtained in terms
of the differentiability properties of $F$.     Proposition \ref{d dimensionale} below is a trivial extension in $d-$variate setting  of Proposition $1.1$ in
Di Bernardino \emph{et al.} (2011)\nocite{DiBernardino}.\\

\begin{Proposition}\label{d dimensionale}
Let  $c \in (0,1)$. Let $F\in\mathcal{F}$ be twice differentiable on $\mathbb{R}^{d*}_{+}$.  Assume there exist $r>0$, $\zeta >0$ such that $m^{\triangledown}>0$ and  $M_H< \infty$.  Then $F$ satisfies  Assumption $\mathbf{H}$,
with $A= \frac{2}{m^{\triangledown}}$.
\end{Proposition}

\vspace{0.1cm}
\begin{remark}\label{curve monotone}
 \rm{Under assumptions of Proposition \ref{d dimensionale}, $\{ F = t \}$ is  a set of Lebesgue measure null in $\mathbb{R}^{d}_{+}$.
Furthermore    we obtain ${\partial L(c)^T= \{ F = c \}^T = \{ F = c \} \cap [0,T]^d}$  (we refer for details to  Remark 1 in Di Bernardino \emph{et al.}, 2011\nocite{DiBernardino} and Theorem 3.2 in Rodr\'{\i}guez-Casal, 2003\nocite{Casal}).}
\end{remark}

\vspace{0.15cm}
\section{Consistency in terms of the Hausdorff  distance}\label{Hausdorff distance}

In this section  we study the consistency properties of $L_n(c)^T$ with respect to the Hausdorff distance between $\partial L_n(c)^T$ and $\partial L(c)^T$. \\

From now on we note, for $n \in \mathbb{N}^*$,
$$\| F -F_n \|_\infty = \sup_{x \, \in \, \mathbb{R}^d_+} \mid F(x)- F_n(x) \mid,$$
and for $T>0$
$$\| F -F_n \|_\infty^T = \sup_{x \, \in \, [0,T]^d} \mid F(x)- F_n(x) \mid.$$

\vspace{0.1cm}
The following result can be considered a trivially  adapted version of Theorem $2.1$ in Di Bernardino \emph{et al.} (2011\nocite{DiBernardino}).

\vspace{0.1cm}
\begin{Theorem}\label{convergenza su insiemi cresenti}
Let  $c \in (0,1)$. Let $F\in\mathcal{F}$ be twice differentiable on $\mathbb{R}^{d*}_{+}$. Assume that there exist $r>0$, $\zeta >0$ such that  $m^{\triangledown}>0$ and $M_H< \infty$. Let $T_1>0$ such that for all
$ \,\, t : \,|\,\, t - c \,\,|\, \leq r,$  $\, \partial L(t)^{T_1} \neq \emptyset$. Let $\left(T_n\right)_{n\in\mathbb{N}^*}$ be an increasing sequence of positive values. Assume that, for each $n$ and for almost all samples of
size $n$, $F_n$ is a continuous function  and that
\begin{equation*}\label{unif convergence bis}
\| F -F_n \|_\infty \rightarrow 0, \quad a.s.
\end{equation*}
Then, for $n$ large enough,
\begin{equation*}\label{Ogrande}
    d_H(\partial L(c)^{T_n}, \partial L_n(c)^{T_n})  \leq 6\, A\, \| F -F_n \|^{T_n}_\infty \,\,\,  \quad  a.s.,
\end{equation*}
where $A=\frac{2}{m^{\triangledown}}$. Therefore we have
\begin{equation*}
d_H(\partial L(c)^{T_n}, \partial L_n(c)^{T_n}) =O(\| F -F_n \|_\infty) \,\,\,  \quad  a.s.
\end{equation*}
\end{Theorem}

\vspace{0.15cm}
Under assumptions of Theorem \ref{convergenza su insiemi cresenti},  $d_H(\partial L(c)^{T_n}, \partial L_n(c)^{T_n})$ converges to zero  and  the quality of our plug-in estimator is obviously related to the quality of the
estimator $F_n$.  For comments and discussions about this result we refer the interested reader  to Remark 2 in Di Bernardino \emph{et al.} (2011\nocite{DiBernardino}).

\section{$L_1$ consistency}\label{lambda consistency}

The previous section was devoted to the consistency of $L_n(c)$ in terms of the Hausdorff  distance. We consider now another consistency criterion: the consistency of the volume (in the Lebesgue measure sense) of the symmetric
difference between $L(c)^{T_n}$ and $L_n(c)^{T_n}$. This means that we define the distance between two subsets $A_1$ and $A_2$ of $\mathbb{R}_d^+$ by

$$d_\lambda(A_1,A_2) = \lambda(A_1\bigtriangleup A_2),$$

\vspace{0.1cm}
where $\lambda$ stands for the Lebesgue measure on $\mathbb{R}^d$ and $\bigtriangleup$ for the symmetric difference.\\

Let us introduce the following assumption:
\begin{itemize}
\item[$\;\mathbf{A1}$]
There exist positive increasing sequences  $\left(v_n\right)_{n\in\mathbb{N}^*}$ and   $\left(T_n\right)_{n\in\mathbb{N}^*}$     such that
$$v_n \, \int_{[0, T_n]^d}\,\mid F -F_n \mid^p \, \lambda( \rmd x ) \mathop{\rightarrow}\limits_{n\to\infty}^\Pbb 0, $$
for some  $1 \leq p < \infty$.
\end{itemize}
\vspace{0.1cm}
We now establish our consistency result with convergence rate, in terms of the volume of the symmetric difference.   We can interpret the following theorem as an extension of Theorem 3.1  in Di Bernardino \emph{et al.}
(2011\nocite{DiBernardino}),   in the case of a $d-$variate distribution function $F$. \\

\vspace{0.1cm}
\begin{Theorem} \label{convergenza lambda su insiemi cresenti}
Let  $c \in (0,1)$. Let $F\in\mathcal{F}$ be  a  twice differentiable distribution function  on $\mathbb{R}^{d*}_{+}$. Assume that there exist $r>0$, $\zeta >0$ such that  $m^{\triangledown}>0$ and $M_H< \infty$. Assume that for
each $n$, with probability one, $F_n$ is measurable.  Let  $\left(v_n\right)_{n\in\mathbb{N}^*}$  and  $\left(T_n\right)_{n\in\mathbb{N}^*}$  positive increasing sequences  such that  Assumption $\mathbf{A1}$  is satisfied and
that for all $ \,\, t : \,|\,\, t - c \,\,|\, \leq r,$   $\, \partial L(t)^{T_1} \neq \emptyset$. Then, it holds that
$$p_n \, d_\lambda(L(c)^{T_n},L_n(c)^{T_n})\mathop{\rightarrow}\limits_{n\to\infty}^{\Pbb}  0,$$
with $p_n$  an  increasing positive sequence such that  $p_n =o\left(v_n^\frac{1}{p+1}/T_n^\frac{(d-1)\,p}{p+1}\right)$.
\end{Theorem}

\vspace{0.1cm}
The proof is postponed to Section \ref{proofs}. This demonstration is basically  based on the proof of Theorem 3.1 in Di Bernardino \emph{et al.}  (2011\nocite{DiBernardino}).\\

Theorem \ref{convergenza lambda su insiemi cresenti} provides a convergence rate, which is closely related to the choice of the sequence $T_n$.
 Note that, as in Theorem 3 in Cuevas \textit{et al.} (2006)\nocite{Cuevas2}, Theorem \ref{convergenza lambda su insiemi cresenti}
above does not require any continuity assumption on $F_n$. Furthermore, as in Theorem 3.1 in Di Bernardino \emph{et al.}  (2011\nocite{DiBernardino}), we remark that a sequence $T_n$, whose divergence rate is large, implies a convergence rate $p_n$ quite slow. Moreover, this phenomenon is emphasized by the dimension $d$ of the data, and we face here the well-known \emph{curse of dimensionality}. In the following we will illustrate this aspect  by giving
convergence rate in the case of the empirical  distribution function (see Example \ref{Example}). Firstly, from Theorem \ref{convergenza lambda su insiemi cresenti} we can derive the following result.

\begin{Corollary}
Under the assumptions and notations of Theorem \ref{convergenza lambda su insiemi cresenti}.  Assume that there exists a positive increasing sequence $(v_n)_{n\in\mathbb{N}^*}$ such that
$v_n \, \| F -F_n \|_\infty  \mathop{\rightarrow}\limits_{n\to\infty}^{\Pbb}  0.$ Then, it holds that
$$p_n \, d_\lambda(L(c)^{T_n},L_n(c)^{T_n})\mathop{\rightarrow}\limits_{n\to\infty}^{\Pbb}  0,$$
with $p_n$  an  increasing positive sequence such that  $p_n =o\left({v_n}^{{\frac {p}{p+1}}}/{T_n}^{{\frac {d+ (d-1)\,p}{p+1}}}\right)$.
\end{Corollary}

This result comes trivially  from Theorem \ref{convergenza lambda su insiemi cresenti} and the fact that that $ v_n \, \| F -F_n \|_\infty  \mathop{\rightarrow}\limits_{n\to\infty}^{\Pbb}  0$ implies

$$\, \forall \,\,\, p \geq 1,  \quad  w_n \, \int_{[0, T_n]^d}\,\mid F -F_n \mid^p \, \lambda( \rmd x ) \mathop{\rightarrow}\limits_{n\to\infty}^\Pbb 0, \quad \mbox{ with  } \quad w_n = \frac{v_{n}^{p}}{T_{n}^{d}}.$$

\vspace{0.1cm}
Let us now present a more practical example.

\vspace{0.1cm}

\begin{Example}[The empirical  distribution function case]\label{Example}$\left.\right.$

Let $F_n$ the $d-$variate empirical distribution function. Then, it holds that  $\, v_n \, \| F -F_n \|_\infty  \mathop{\rightarrow}\limits_{n\to\infty}^{\Pbb}  0,\,\,$ with  $v_n=o(\sqrt{n})$. From Theorem
\ref{convergenza lambda su insiemi cresenti}, with $p=2$,  we obtain   for instance:

 $$ p_n =o\left({\frac {{n}^{1/3}}{{T_n}^{7/3}}}\right), \quad  \mbox{ for } \,\, d=3; \quad \quad \quad
p_n =o\left({\frac {{n}^{1/3}}{{T_n}^{10/3}}}\right), \quad  \mbox{ for }  \,\, d= 4.$$
\end{Example}

$\mbox{}$
\vspace{0.1cm}

 The next section is dedicated to study the effects of scaling  data.

\section{About the effects of scaling data}\label{variation echelle}

Suppose now to scale our data using a scale parameter $a \in \Real^*_+$. In our case, the scaled random vector will be  $(a \,X_1,a \,X_2,\ldots,a \,X_d):=a \, \textbf{X}$. From now on we denote $F_{a \,\textbf{X}}$
(resp.  $F_{\textbf{X}}$) the distribution  function associated to $a \, \textbf{X}$ (resp. to $\textbf{X}$).  Using  notation of Section \ref{Notations}, let
$$L_{a}(c) = \{ x \in \Real^d_{+}  : F_{a \, \textbf{X}}(x) \geq c \}.$$

It is easy to prove (see for instance Section 3 in Tibiletti, 1993\nocite{Tibiletti3}) that
$$L_{a}(c)= a \, L(c),$$
and
$$E_{a} = B(\{x  \in \mathbb{R}^{d}_{+} :\,  \mid F_{a \, \textbf{X}}(x)-c\mid \leq r\}, \, \zeta)= a\,E.$$

Define now
\begin{equation*}
m^{\triangledown}_{a} =\inf_{x  \in  E_a} \| \nabla F_{a \,\textbf{X}}(x)\|.
\end{equation*}

First, we can obtain the following result whose proof is postponed to Section \ref{proofs}.

\begin{Lemma}\label{m deltaR}
It holds that
  \begin{equation*}\label{nuovo mdeltaR}
  m^{\triangledown}_{a}= \frac{1}{a}\, m^{\triangledown}, \quad \forall \,\, a \in \Real^*_+.
  \end{equation*}
  Furthermore, if
$$ M_H = \sup_{x \in E} \|(H  F_{\textbf{X}})_{x}\| < + \infty  \quad  \mbox{ then } \quad   M_{H, a} = \sup_{x \in a\,E} \|(H  F_{a \,\textbf{X}})_{x}\|< + \infty,  \quad \mbox{ with } \,\, a \in \Real^*_+.$$
\end{Lemma}

  \vspace{0.1cm}

We can now consider the effects  of scaling  data on Theorem \ref{convergenza su insiemi cresenti} and \ref{convergenza lambda su insiemi cresenti}.

\begin{Theorem}\label{meme A}

$\mbox{ }$

\begin{enumerate}
\item Under same notation and assumptions of Theorem \ref{convergenza su insiemi cresenti}, for $n$ large enough, it holds that
\begin{equation*}\label{Ogrande}
    d_H(\partial L_{a}(c)^{a\,T_n}, \partial L_{n,\,a}(c)^{a\,T_n})  \leq 6\, A\, a \, \| F -F_n \|^{T_n}_\infty, \,\,\,  \quad  a.s.
\end{equation*}
  \vspace{0.12cm}
\item Under same notation and assumptions of Theorem \ref{convergenza lambda su insiemi cresenti} it holds that
 $$p_{n, \,a} \, d_\lambda(L_a(c)^{a\,T_n},L_{n,\, a}(c)^{a\, T_n})\mathop{\rightarrow}\limits_{n\to\infty}^{\Pbb}  0,$$
with $p_{n, \,a}$  an  increasing positive sequence such that  $p_{n, \,a} =o\left(v_n^\frac{1}{p+1}/\left(a^\frac{d\,p}{p+1} \, T_n^\frac{(d-1)\,p}{p+1}\right)\right)$.
\end{enumerate}
\end{Theorem}

\begin{remark}\label{pna}

$\mbox{ }$

\begin{enumerate}
{\rm
\item The first result of Theorem \ref{meme A} states that a change of scale of the data implies the same change of scale for the Hausdorff distance.\\

\item  The second result states that a change of scale  of the data implies a rate in $$o\left(v_n^\frac{1}{p+1}/\left(a^d\,T_n^{(d-1)}\right)^{p/(p+1)}\right)$$
instead of $$o\left(v_n^\frac{1}{p+1}/ \left(T_n^{(d-1)}\right)^{\frac{p}{p+1}}\right).$$
So, we see logically that the scale factor $a$ impacts the volume in $\mathbb{R}^d$ with an exponent $d$.
}
\end{enumerate}
\end{remark}

\section*{Conclusion}\label{Conclusions}

Starting from previous results obtained in Di Bernardino \emph{et al.} (2011\nocite{DiBernardino}), we propose in this paper a generalization to the estimation of level sets in the case of a $d$-variate distribution function. The consistency
results are stated in term of Hausdorff distance and volume of the symmetric difference. We propose a rate of convergence for this second criterion. Moreover, we analyze the impact of scaling  data  on our results. As a
future work, a complete simulation study and an \verb"R"-package are in preparation.

\section{Proofs}\label{proofs}

\textbf{\small Proof of Theorem \ref{convergenza lambda su insiemi cresenti}\vspace{0.2cm}}

Under assumptions of  Theorem \ref{convergenza lambda su insiemi cresenti},  we can always take  $T_1 > 0$ such that  for all $ \,\,t : \,|\, t - c \, |\leq r,$  $\, \partial L(t)^{T_1} \neq \emptyset$.  Then for each $n$,  for all
$ \,\,t : |\, t - c \, |\leq r,$    $\, \partial L(t)^{T_n}$ is a non-empty  (and compact) set on $\mathbb{R}^d_+$.\\

\begin{sloppypar}
We consider a positive sequence $\varepsilon_n$ such that  $\varepsilon_n \mathop{\rightarrow}\limits_{n\to\infty} 0$. For each $n\geq1$ the random sets
${L(c)^{T_n} \bigtriangleup L_n(c)^{T_n}}$, ${Q_{\varepsilon_n} = \{x \in   [0,T_n]^d :\, \mid F -F_n \mid \leq \varepsilon_n\}}$ and ${\widetilde{Q}_{\varepsilon_ n} = \{x \in   [0,T_n]^d: \,\, \mid F -F_n \mid >\varepsilon_n\}}$
are measurable and
\begin{equation*}
\lambda(L(c)^{T_n} \bigtriangleup   L_n(c)^{T_n})= \lambda(L(c)^{T_n} \bigtriangleup  L_n(c)^{T_n}\cap \, Q_{\varepsilon_n})+ \lambda(L(c)^{T_n} \bigtriangleup  L_n(c)^{T_n}\cap \,{\widetilde{Q}}_{\varepsilon_n}).
\end{equation*}
\end{sloppypar}
Since $L(c)^{T_n} \bigtriangleup L_n(c)^{T_n}\cap \, Q_{\varepsilon_n}\subset \{x \in  [0,T_n]^d : c- \varepsilon_n \leq F <   c + \varepsilon_n\}$ we obtain
$$\lambda(L(c)^{T_n} \bigtriangleup  L_n(c)^{T_n}) \leq \lambda(\{x \in [0, T_n]^d: c- \varepsilon_n \leq F <   c + \varepsilon_n\})+ \lambda({\widetilde{Q}}_{\varepsilon_n}).$$

From Assumption $\mathbf{H}$ (Section \ref{Notations}) and Proposition \ref{d dimensionale}, if $2 \,\varepsilon_n \leq \gamma$ then \begin{center}$d_H(\partial L(c+\varepsilon_n)^{T_n},\partial L(c-\varepsilon_n)^{T_n}) \leq  2 \, \varepsilon_n \, A.$\end{center}

From assumptions on first derivatives of $F$ (see Assumption $\mathbf{H}$  and Proposition \ref{d dimensionale}) and Propriety 1 in  Imlahi \emph{et al.} (1999)\nocite{Imlahi},   we can write
$$\lambda(\{x \in [0, T_n]^d: c- \varepsilon_n \leq F <   c + \varepsilon_n\}) \leq \,  (2 \, \varepsilon_n \, A) \, d \, T_n^{\,d-1}.$$

If we now choose
\begin{equation} \label{epsilon1}
 \varepsilon_n=o\left(\frac{1}{p_n \, T_n^{\,d-1}}\right),
\end{equation}
we obtain that, for $n$ large enough,   $2 \,\varepsilon_n \leq \gamma$ and
$$p_n \, \lambda(\{x \in  [0,T_n]^d : c- \varepsilon_n \leq F <   c + \varepsilon_n\})\, \mathop{\rightarrow}\limits_{n\to\infty} \, 0. $$

Let us now prove that  $p_n \, \lambda({\widetilde{Q}}_{\varepsilon_ n}) \, \mathop{\rightarrow}\limits_{n\to\infty}^{\Pbb} \, 0.$ To this end, we write
\begin{equation*}
 p_n \, \lambda({\widetilde{Q}}_{\varepsilon_ n}) \,= \,  p_n  \int \, 1_{\{x \in  [0,T_n]^d : \, \mid F -F_n \mid > \varepsilon_n\}} \, \lambda(\rmd x ) \leq \, \frac{p_n}{\varepsilon_n^p}  \, \int_{[0, T_n]^d}   \mid F -F_n \mid^p \, \lambda(\rmd x ).
\end{equation*}
Take $\varepsilon_n$ such that  \begin{equation}\label{calcolo stimatore generico}\varepsilon_n= \, \left({\frac{p_n}{v_n}}\right)^\frac{1}{p}.\end{equation}
So, from Assumption $\mathbf{A1}$ in Section \ref{lambda consistency}, we obtain $p_n \, \lambda({\widetilde{Q}}_{\varepsilon_ n}) \, \mathop{\rightarrow}\limits_{n\to\infty}^{\Pbb} \,0.$  As  $p_n =o\left(v_n^\frac{1}{p+1}/T_n^\frac{(d-1)\,p}{p+1}\right)$  we can choose $\varepsilon_n$ that satisfies  \eqref{epsilon1} and   \eqref{calcolo stimatore generico}.   Hence the result.  $\, \square$ \vspace{0.3cm}

\textbf{\small Proof of Lemma \ref{m deltaR}\vspace{0.2cm}}

First, we remark that
$$ F_{a \, \textbf{X}}(x)= F_{\textbf{X}}\left(\frac{x}{a}\right).$$
Then, we obtain 

\begin{eqnarray*}
    m^{\triangledown}_{a}& = &\inf_{x  \in a\,E} \left\| \left( \frac{\partial}{\partial x_1}F_{\textbf{X}}\left(\frac{x}{a}\right), \ldots, \frac{\partial}{\partial x_d}F_{\textbf{X}}\left(\frac{x}{a}\right) \right) \right\|, \\
  &=& \inf_{x  \in a\,E} \left\| \frac{1}{a} \left(\frac{\partial F_{\textbf{X}}}{\partial x_1}\left(\frac{x}{a}\right), \ldots, \frac{\partial F_{\textbf{X}}}{\partial x_d}\left(\frac{x}{a}\right)\right) \right\|, \\
  &= & \frac{1}{a} \, \inf_{x  \in  \,E} \left\|   \left(\frac{\partial F_{\textbf{X}}}{\partial x_1}(x), \ldots, \frac{\partial F_{\textbf{X}}}{\partial x_d}(x)\right) \right\|.  \\
  &= & \frac{1}{a}\,  m^{\triangledown}.
\end{eqnarray*}
Second part of Lemma  \ref{m deltaR} comes down from trivial calculus. Hence the result.  $\, \square$ \vspace{0.3cm}

\textbf{\small Proof of Theorem  \ref{meme A}\vspace{0.2cm}}

\textit{$\qquad$ Proof of $1$.}\vspace{0.1cm}

 Following the proof of Theorem \ref{convergenza su insiemi cresenti},  it holds that

 $$d_H(\partial L_{a \, \textbf{X}}(c)^{a\,T_n}, \partial L_{n,\,a}(c)^{a\,T_n})  \leq 6\, \frac{2}{m^{\triangledown}_{a}} \, \,\, \sup_{x \in [0, a\,Tn]^d}
   \mid   F_{\textbf{X}}\left(\frac{x}{a}\right) -F_{n}\left(\frac{x}{a}\right) \mid.$$

Using Lemma \ref{m deltaR} and the fact that
$$  \sup_{x \in [0, a\,Tn]^d}
   \mid   F_{\textbf{X}}\left(\frac{x}{a}\right) -F_{n}\left(\frac{x}{a}\right) \mid = \sup_{x \in [0, Tn]^d}
   \mid   F_{\textbf{X}}(x) -F_{n}(x) \mid, $$
 we get the result. $\, \square$ \vspace{0.3cm}

\textit{$\qquad$ Proof of $2$.}\vspace{0.1cm}

As in the proof of Theorem \ref{convergenza lambda su insiemi cresenti} and using  same notation, we can write
$$\lambda(\{x \in [0, a\, T_n]^d: c- \varepsilon_n \leq F_{a\, \textbf{X}} <   c + \varepsilon_n\}) \leq \,  (2 \, \varepsilon_n \, A \,a) \,  d \, a^{d-1}\,T_n^{d-1}.$$

 If we now choose
 \begin{equation} \label{epsilon primo}
 \varepsilon_n=o\left(\frac{1}{p_{n, \,a} \, a^{d}\, T_n^{d-1}}\right)\end{equation}
we obtain that for $n$ large enough   $2 \,\varepsilon_n \leq \gamma$ and
$$p_{n, \,a} \, \lambda(\{x \in [0, a\, T_n]^d: c- \varepsilon_n \leq F_{a\, \textbf{X}} <   c + \varepsilon_n\})  \mathop{\rightarrow}\limits_{n\to\infty} 0. $$

The second part of this demonstration is equal to proof of Theorem \ref{convergenza lambda su insiemi cresenti}.  Then we take
$\varepsilon_n$ such that  \begin{equation}\label{epsilon secondo}\varepsilon_n=  \left({\frac{p_{n, \,a}}{v_n}}\right)^\frac{1}{p}.\end{equation}  Then,  from Assumption $\mathbf{A1}$, in Section \ref{lambda consistency}, we obtain   $p_{n, \,a} \, \lambda(\{x \in   [0,a\,T_n]^d: \,\, \mid F_{a\,\textbf{X}} -F_{a,\,n} \mid >\varepsilon_n\}) \, \mathop{\rightarrow}\limits_{n\to\infty}^{\Pbb} 0.$  As $p_{n, \,a} =o\left(v_n^\frac{1}{p+1}/\left(a^\frac{d\,p}{p+1} \, T_n^\frac{(d-1)\,p}{p+1}\right)\right)$ we can choose $\varepsilon_n$ that satisfies  \eqref{epsilon primo} and   \eqref{epsilon secondo}.   Hence the result.  $\, \square$ \\\\

\vspace{0.15cm}
\noindent
\textbf{Acknowledgements: }  This work has been partially supported by the French research national agency (ANR) under the reference ANR-08BLAN-0314-01. The authors thank  Yannick Baraud, Christine Tuleau-Malot and  Patricia Reynaud-Bouret   for fruitful discussions. \\\\

\vspace{0.15cm}
\noindent
\textbf{References: }
 \bibliographystyle{abbrv}\footnotesize
\addcontentsline{toc}{section}{Bibliography}
\bibliography{biblio}

\end{document}